\documentclass[11pt]{article}
\usepackage{amsmath,amsmath,amsthm,amssymb}
\usepackage{courier}
\usepackage{graphicx}
\usepackage{esint}
\usepackage{enumerate}
\usepackage{fullpage}
\usepackage{hyperref}

\newtheorem{theorem}{Theorem}
\newtheorem{lemma}{Lemma}
\newtheorem{corollary}{Corollary}

\numberwithin{equation}{section}

\begin{document}

\title{On problems with weighted elliptic operator and general growth nonlinearities}

\author{John Villavert\footnote{email: john.villavert@gmail.com, john.villavert@utrgv.edu} \\
[0.2cm] {\small University of Texas, Rio Grande Valley}\\
{\small Edinburg, TX 78539, USA}
}

\date{}
 

\maketitle

\begin{abstract} 
This article establishes existence, non-existence and Liouville-type theorems for nonlinear equations of the form
\begin{equation*}
-div (|x|^{a} D u ) = f(x,u), ~ u > 0,\, \mbox{ in } \Omega,
\end{equation*}
where $N \geq 3$, $\Omega$ is an open domain in $\mathbb{R}^N$ containing the origin, $N-2+a > 0$ and $f$ satisfies structural conditions, including certain growth properties. The first main result is a non-existence theorem for boundary-value problems in bounded domains star-shaped with respect to the origin, provided $f$ exhibits supercritical growth. A consequence of this is the existence of positive entire solutions to the equation for $f$ exhibiting the same growth. A Liouville-type theorem is then established, which asserts no positive solution of the equation in $\Omega = \mathbb{R}^N$ exists provided the growth of $f$ is subcritical. The results are then extended to systems of the form
\begin{equation*}
-div (|x|^{a} D u_1) \!=\! f_{1}(x,u_1,u_2),  -div (|x|^{a} D u_2) \!=\! f_{2}(x,u_1,u_2),  u_1, u_2 \!>\! 0,\, \mbox{ in } \Omega,
\end{equation*}
but after overcoming additional obstacles not present in the single equation. Specific cases of our results recover classical ones for a renowned problem connected with finding best constants in Hardy-Sobolev and Caffarelli-Kohn-Nirenberg inequalities as well as existence results for well-known elliptic systems.
\end{abstract}

 \medskip


\section{Introduction}\label{Introduction}

In this paper, we first examine the general elliptic equation
\begin{equation}\label{pde}
-div (|x|^{a} D u ) = f(x,u), ~ u > 0,\, \mbox{ in } \Omega,
\end{equation}
where $N \geq 3$, $\Omega$ is an open domain containing the origin, $\Omega_0 := \Omega \backslash\{0\}$, $p > 1$, $N - 2 + a > 0$, $Du$ denotes the gradient of $u$, and $f: \Omega_0 \times [0,\infty) \longrightarrow \mathbb{R}$ is a smooth mapping satisfying a local Lipschitz-type condition and various structural and growth assumptions. Specifically, we always assume that
\begin{enumerate}[(I)]
\item $f(x,0) = 0$ and $f(x,u) > 0$ for all $x \in \Omega_0$, $u > 0$,
\item $0 \leq u \leq v$ implies $f(x,u) \leq f(x,v)$, for all $x \in \Omega_0$,
\item there is a $b > a - 2 > -N$ such that for each compact $\omega \subset [0,\infty)$, there exists $L>0$ such that
\[ |f(x,u) - f(x,v)| \leq L|x|^b |u - v| ~ \mbox{ for } x \in \Omega_0,\, u, v \in \omega.\]
\end{enumerate}
The methods we employ for the single equation, after careful modifications and under certain assumptions, do carry over to systems of the form
\begin{align}\label{pde system}
         &div(|x|^a D u_1) \!+\! f_{1}(x,u_1,u_2) \!=\! 0, \nonumber \\ &div(|x|^a D u_2) \!+\! f_{2}(x,u_1,u_2) \!=\! 0,  u_{1}, u_{2} \!>\! 0,  \mbox{ in } \Omega.
\end{align}
So for completeness, we shall include results for this system akin to those for single equations. For simplicity we will often express \eqref{pde system} in the concise vector form
\begin{equation}\label{vector system}
 div(|x|^a DU) + F(x,U) = 0, \, U > 0, \,\mbox{ in } \Omega,
\end{equation}
where $U = (u_1,u_2)$ and $F = (f_1,f_2)$. Relations such as $U > 0$, $U \geq 0$, etc, are understood to hold component-wise. If $U > 0$ and $P = (p_1,p_2) \geq (1,1)$, we define the scalars $|P| = p_1 + p_2$ and $U^P = u_{1}^{p_1}u_{2}^{p_2}$. We always assume $f_1, f_2: \Omega_0 \times [0,\infty)\times [0,\infty) \longrightarrow [0,\infty)$ are both smooth and satisfy the analogues of conditions (I)--(III) as above for $f$, except we replace the second condition with the {\bf cooperative property}: for $i = 1,2$,
\begin{equation}\label{cooperative}
\text{(II')} \qquad  U \leq V \mbox{ and } u_i = v_i \mbox{ imply } f_{i}(x,U) \leq f_{i}(x,V).
\end{equation}
By a positive solution $u$ of problem \eqref{pde} we mean a function $u:\Omega \longrightarrow (0,\infty)$ of the class $C^{2}(\Omega_0) \cap C(\overline{\Omega})$ satisfying the equation pointwise everywhere in $\Omega_0 := \Omega \backslash \{0\}$. A solution $U = (u_1,u_2) \in (C^{2}(\Omega_0) \cap C(\overline{\Omega}) )^2$ of system \eqref{pde system} is understood analogously.

Let us motivate our interest in studying equation \eqref{pde} and system \eqref{pde system}. For the scalar equation, the model example here is perhaps the weighted power nonlinearity $f(x,u) = |x|^b u^p$, which reduces \eqref{pde} into the well-known problem,
\begin{equation}\label{ckn}
div( |x|^{a} D u) + |x|^{b}u^p = 0, ~ u> 0, \mbox{ in } \Omega,
\end{equation}
where $N \geq 3$, $p > 1$, $b > - N$ and $N - 2 + a > 0$. Equation \eqref{ckn} is essentially the Euler-Lagrange equation for a variational problem connected with the following sharp Caffarelli-Kohn-Nirenberg inequality \cite{CKN84}: there exists a best constant $C = C(N,a,b)> 0$, depending only on $N$, $a$ and $b$, such that
\begin{equation*}
C \Big( \int_{\mathbb{R}^N} |x|^b |f|^{p+1} \,dx\Big)^{2/(p+1)} \leq \int_{\mathbb{R}^N} |x|^a |Df|^{2} \,dx \mbox{ for } f \in C_{c}^{\infty}(\mathbb{R}^N),
\end{equation*}
provided that $\frac{N+b}{p+1} + 1 = \frac{N+a}{2}$. Quantitative and qualitative properties--e.g., existence, non-existence, classification, asymptotic behavior, symmetry-breaking, etc--of the optimizers for this class of sharp inequalities were examined earlier in \cite{CW01,ChouChu93} (also see the cited papers there). Additional examination of the differential equation itself can be found in \cite{DuGuo13,DuGuo15,GuoWan17,Villavert:20} and the references therein. Our results may be viewed as direct extensions of some of those obtained in these papers. Further, the motivating problems mentioned above and studied in the cited papers also explain some of the natural assumptions we place on the parameters $N$, $a$, $b$, $p$, etc.

We should mention the notable case $a = 0$, which gives the equation
\begin{equation}\label{hle}
          \Delta u + |x|^{b} u^p = 0, \, u > 0, \, \mbox{ in } \Omega.
\end{equation}
Problem \eqref{hle} is sometimes known as the H\'{e}non-Lane-Emden equation (or just Lane-Emden equation if $b=0$). If $\Omega = \mathbb{R}^N$, equation \eqref{hle} has several important applications, e.g., it arises as an astrophysical model for stellar cluster formation; it comprises the blow-up equation used to obtain a priori estimates for a class of elliptic boundary value problems; and it appears in geometric problems such as Yamabe's problem and sharp Sobolev inequalities. For additional closely-related papers concerning equation \eqref{ckn}, the reader is referred to \cite{BVG10,DDG11,DuGuo13} and the references therein.

For system \eqref{pde system}, specific cases such as those involving Laplace operators and weighted power nonlinearities have received considerable attention in the past few decades. Much effort has been directed in obtaining sharp existence and non-existence results on the positive solutions \cite{BM02,CL09,FG14,Phan12,RZ00}. Obtaining such sharp existence results is typically far more difficult for coupled systems than for scalar equations. For instance, we mention the renowned Lane-Emden system, which corresponds to \eqref{pde system} when $N\geq 2$, $a=b= 0$, $p_1, \, p_2 > 0$, $f_{1}(x,U) = u_{2}^{p_1}$ and $f_{2}(x,U) = u_{1}^{p_2}$. Even for this relatively simpler case, establishing an optimal Liouville-type theorem remains open though partial results are known. More precisely, the conjecture asserts the Lane-Emden system admits no positive classical solutions if and only if
\begin{equation}\label{subcritical hyperbola}
H(p_1,p_2) := \frac{1}{1+p_{1}} + \frac{1}{1+p_{2}} > \frac{N-2}{N}.
\end{equation}
The existence of positive radial solutions whenever $H(p_1,p_2) \leq (N-2)/N$ was obtained in \cite{SZ98}. Meanwhile, the non-existence of positive radial solutions whenever \eqref{subcritical hyperbola} holds was proved in \cite{Mitidieri96}. Thus, the conjecture is resolved in the class of radially symmetric solutions. In general, however, the conjecture--specifically the non-existence part--has only been verified in the lower dimensions $N\leq 4$ \cite{PQS07,SZ96,Souplet09} (see also \cite{BM02,DeFF94}).

Our main motivation in this article is to identify conditions ensuring the existence or non-existence of positive solutions to problem \eqref{pde} for both $\Omega = \mathbb{R}^N$ and for bounded star-shaped domains with prescribed zero boundary conditions. Consequently, we recover and extend such results for problem \eqref{ckn} and its special cases. As our results illustrate, albeit not surprisingly, the existence and non-existence of positive solutions depend precisely on the domain and the structure and growth behavior of the nonlinearities. Furthermore, an underlying difficulty in examining our class of problems stems from the possible singular or degenerate operator $div(|x|^a D\, \cdot)$ coupled with a possible inhomogeneity in the nonlinearity $f$, see, e.g., equation \eqref{ckn}. The previous methods which proved successful for addressing specific cases, e.g. the Lane-Emden equation, no longer apply easily to our general problem. We circumvent such difficulties in obtaining non-existence results using ideas influenced mainly by the work in \cite{BM02,RZ00}, and for the existence results we adopt ideas from \cite{LV16a,LV16b,LGZ06a,Villavert:14b}. That is, by carefully utilizing shooting methods, Kelvin transforms, comparison arguments, and moving plane methods adapted to infinite cylindrical domains, we are able to extend the classical existence and non-existence results to equation \eqref{pde} and system \eqref{pde system}.

\subsection{Results for single equations}

Our first main result is a non-existence theorem for bounded domains, provided that $f$ has supercritical growth and $\Omega$ is star-shaped.

\begin{theorem}\label{thm1}
Let $N \geq 3$, $N - 2 + a > 0$, $\Omega$ is an open bounded domain star-shaped with respect to the origin, and $\partial \Omega$ is smooth. We further assume for every $(x,y) \in \Omega_0 \times [0,\infty)$,
\begin{equation}\label{supercritical}
\mu \longrightarrow \mu^{-\frac{N+2-a}{N-2+a}}f(\mu^{-\frac{2}{N-2+a}}x, \mu y) ~\mbox{ is non-decreasing in } \, [1,\infty).
\end{equation}
Then the problem
\begin{equation}\label{gen bvp}
  \begin{cases}
          div( |x|^{a} D u) + f(x,u) = 0  \quad \text{ in } \Omega, \\
 	  u= 0  \qquad \qquad \qquad \qquad \quad \;\;\text{ on } \partial\Omega,
        \end{cases}
\end{equation}
has no positive solution in $\Omega$.
\end{theorem}
The above theorem is reminiscent of the non-existence result for the Lane-Emden equation, however, we prove it in the spirit of \cite{RZ00} via comparison arguments and Kelvin transforms rather than through Pohozaev-type identities. Combining Theorem \ref{thm1} with a basic shooting argument yields the following existence result for \eqref{pde} in $\Omega = \mathbb{R}^N$.
\begin{theorem}\label{thm2}
Let $\Omega = \mathbb{R}^N$, $N \geq 3$ and $N - 2 + a > 0$. We further assume for every $(x,y) \in \mathbb{R}^N \backslash \{0\} \times [0,\infty)$, the growth condition \eqref{supercritical} holds. Then there exists a positive solution to equation \eqref{pde}.
\end{theorem}
For the whole space domain, we have the following Liouville-type theorem for problem \eqref{pde}.

\begin{theorem}\label{thm3}
Let $\Omega = \mathbb{R}^N$, $N \geq 3$, $N-2+a > 0$ and the nonlinearity $f$ satisfies
\begin{enumerate}[(i)]
\item for every $(x,y) \in \Omega_0 \times [0,\infty)$,
\begin{equation}\label{subcritical}
\mu \longrightarrow \mu^{\frac{N+2-a}{N-2+a}}f(\mu^{\frac{2}{N-2+a}}x, \mu^{-1}y) ~\mbox{ is (strictly) increasing in } \, \mu \in [1,\infty),
\end{equation}
\item there exist $d > 0$, $b > a - 2$ and $1 < p < p_{S}(a,b):= \frac{N+2+2b-a}{N-2+a}$ such that for large $R>1$,
\begin{equation}\label{supersoln}
 f(x,y) \geq d |x|^b y^p ~ \mbox{ for } (x, y) \in B_{R}^{c}(0) \times [0,\infty).
 \end{equation}
\end{enumerate}
Then problem \eqref{pde} has no positive solution.
\end{theorem}

Nonlinearities exhibiting properties \eqref{supercritical} and \eqref{subcritical} are sometimes said to have supercritical and (strictly) subcritical growth, respectively. Our general results above recover and extend the following sharp existence result for equation \eqref{ckn} with $\Omega = \mathbb{R}^N$.

\begin{corollary}\label{cor1}
Let $\Omega = \mathbb{R}^N$, $N \geq 3$, $p > 1$ and $b > a - 2 > -N$.
\begin{enumerate}[(a)]
\item Equation \eqref{ckn} has no positive solution if $ p < p_{S}(a,b).$
\item Equation \eqref{ckn} admits a positive solution if $p \geq p_{S}(a,b)$.
\end{enumerate}
\end{corollary}
We note that the existence result of (a) and the non-existence result of (b) in Corollary \ref{cor1} were already obtained, among other interesting and related results, in \cite{GuoWan17} and \cite{DuGuo13}, respectively.

\subsection{Results for the system of equations}

We extend the previous results to system \eqref{pde system}. Although we establish these results via similar approaches from the scalar case, obstructions appear in the systems case requiring us to place additional assumptions on nonlinearities and to modify some of the methods we employ.

\begin{theorem}\label{thm4}
Let $N \geq 3$, $N - 2 + a > 0$, $E = [0,\infty) \times [0,\infty)$ and assume for every $(x,Y) \in \Omega_0 \times E$,
\begin{equation}\label{supercritical sys}
\mu \longrightarrow \mu^{-\frac{N+2-a}{N-2+a}}F(\mu^{-\frac{2}{N-2+a}}x, \mu Y) ~\mbox{ is non-decreasing in } \, [1,\infty).
\end{equation}

(a) Let $\Omega$ be an open, bounded smooth domain that is star-shaped with respect to the origin. Then the problem
\begin{equation}\label{gen sys bvp}
 \begin{cases}
          div( |x|^{a} D U) + F(x,U) = 0  \quad \text{ in } \Omega, \\
 	  U = 0 \qquad \qquad \qquad \qquad \qquad  \text{ on } \partial\Omega, 
        \end{cases}
\end{equation}
has no positive solution in $\Omega$.
(b) If $\Omega = \mathbb{R}^N$ and it further holds that
\begin{equation}\label{source compare}
u_{i} \leq u_{j} \, (i\neq j) ~\mbox{ implies } ~ f_{j}(x,U) \leq f_{i}(x,U),
\end{equation}
then \eqref{vector system} admits a positive solution.
\end{theorem}

\begin{theorem}\label{thm5}
Let $\Omega = \mathbb{R}^N$, $N \geq 3$, $N-2+a > 0$ and the nonlinearity $F$ satisfies

(i)for every $(x,Y) \in \Omega_0 \times E$,
\begin{equation}\label{subcritical sys}
\mu \longrightarrow \mu^{\frac{N+2-a}{N-2+a}}F(\mu^{\frac{2}{N-2+a}}x, \mu^{-1}Y) ~\mbox{ is (strictly) increasing in } \, \mu \in [1,\infty).
\end{equation}

(ii) For $i =1,2$ there exist $d  > 0$, $b > a - 2$ and $P^{i} > (1,1)$ with $|P^i | < p_{S}(a,b)$ such that for large $R>1$,
\begin{equation}\label{supersoln sys}
 f_{i}(x,Y) \geq d |x|^b Y^{P^i} ~ \mbox{ for } (x, Y) \in B_{R}^{c}(0) \times E.
 \end{equation}
Then problem \eqref{vector system} has no positive entire solution.
\end{theorem}

Consequences of Theorems \ref{thm4} and \ref{thm5} are the following.
\begin{corollary}
Let $N \geq 3$, $\Omega \subseteq \mathbb{R}^N$, $p,q,r,s > 1$, $b > a - 2 > -N$, and consider the system
\begin{equation}\label{model system}
 \begin{cases}
         div(|x|^a D u_1) + |x|^b u_{1}^{p}u_{2}^{q} = 0,  \\
         div(|x|^a D u_2) + |x|^b u_{1}^{r}u_{2}^{s}  = 0,  u_1, u_{2} > 0,  \mbox{ in } \Omega. 
        \end{cases}
\end{equation}

(a) For $\Omega = \mathbb{R}^N$, system \eqref{model system} has no positive solution whenever $$\max\{p + q, r + s\} < p_{S}(a,b).$$

(b) Let $\Omega$ be a bounded, open smooth domain star-shaped with respect to the origin. Then system \eqref{model system} with the boundary conditions, $u_1 = u_2 = 0$ on $\partial \Omega$, has no positive solution in $\Omega$ whenever
\begin{equation*}
\min\{p + q, r + s\} \geq p_{S}(a,b).
\end{equation*}
For $\Omega = \mathbb{R}^N$, system \eqref{model system} admits a positive solution if $q - s = r - p \geq 1$ and $r + s \geq p_{S}(a,b)$.
\end{corollary}

This paper is organized as follows. In Section \ref{preparations}, we arrive at several intermediate results that comprise the essential ingredients in our proofs of the Liouville-type theorems. Particularly, after rewriting the equation via an Emden-Fowler type transformation, we then apply a moving planes approach to get a key monotonicity result. In Section \ref{penultimate section}, we provide the proof of Theorems \ref{thm1}--\ref{thm2} and Theorem \ref{thm4}, where the existence results in these theorems will follow from a basic shooting argument combined with our non-existence result on bounded star-shaped domains. We should remark that this shooting method requires somewhat of a technical refinement when adapted for the system case. Then Section \ref{final section} contains the proofs of Theorems \ref{thm3} and Theorem \ref{thm5}, which relies closely on the monotonicity result and a comparison with rescaled sub-solutions to a model boundary-value problem.


\section{Preparations and intermediate results}\label{preparations}

Let us first discuss the notation and conventions we adopt hereafter. We denote by $B_{R}(x) \subset \mathbb{R}^N$ the open ball of radius $R>0$ centered at $x \in \mathbb{R}^N$. We denote its boundary by $\partial B_{R}(x)$, and if $x = 0$ and $R = 1$, then we write the resulting $(N-1)$-dimensional unit sphere $\partial B_{1}(0)$ as $\mathbb{S}^{N-1}$ instead. We denote the complement of $B_{R}(x)$ by $B_{R}^{c}(x) = \mathbb{R}^N \backslash B_{R}(x)$. The constant $C$ in inequalities below represents some universal constant that may change from line to line, or even within the same line itself.

Some of our methods will occasionally depend on writing \eqref{pde} and its related equations in polar coordinates. Namely, if $u$ is a solution of \eqref{pde} in $\Omega = \mathbb{R}^N$ and, for every non-zero $x$, we write
\begin{equation}\label{polar coordinates}
r = |x|  \mbox{ and }  \theta = x/|x| \in \mathbb{S}^{N-1},
\end{equation}
and $ u(x) = v(r,\theta)$. Then
\begin{equation*}
 div(|x|^a Du) = r^{a}\Big( \partial_{r}^2 v + \frac{N - 1 + a}{r}\partial_{r}v + \frac{1}{r^2}\Delta_{\theta}v\Big),
\end{equation*}
where $\partial_{r}^k := \frac{\partial^k}{\partial r^k}$ and $\Delta_{\theta}$ is the Laplace-Beltrami operator on $\mathbb{S}^{N-1}$. Writing $f(x,u)$ as $f(r,\theta,v)$, it follows that $v = v(r,\theta)$ satisfies
\begin{equation}\label{polar equation}
r^a( \partial_{r}^2 v + \frac{N - 1 + a}{r}\partial_{r}v + \frac{1}{r^2}\Delta_{\theta}v) + f(r,\theta,v ) = 0, ~ v > 0, ~ \mbox{ in } (0,\infty) \times \mathbb{S}^{N-1}.
\end{equation}
In view of this, we will make use of the notation
$$L[u] = \Delta u + a|x|^{-2}(x\cdot Du) \, \mbox{ and } \, L_{r}[v] = \partial_{r}^2 v + \frac{N - 1 + a}{r}\partial_{r}v. $$


\subsection{A monotonicity property}
A key ingredient to establishing the Liouville-type theorems is the following monotonicity result, which we derive via the adapted method of moving planes.
\begin{lemma}\label{lemma4}
Let $N \geq 3$, $\Omega = \mathbb{R}^N$, $\gamma = (N-2+a)/2 > 0$, and suppose $u$ is a positive solution of \eqref{pde} and that assumption (i) in Theorem \ref{thm3} holds. Then $|x|^{\gamma}u(x)$ is monotone increasing with respect to $|x|$.
\end{lemma}

\begin{proof}
Let $u$ be a positive solution of \eqref{pde} with $\Omega = \mathbb{R}^N$.

{\bf Step 1.} We apply an Emden-Fowler type transformation.

Let $(r,\theta) \in (0,\infty) \times \mathbb{S}^{N-1}$ represent polar coordinates as defined in \eqref{polar coordinates} and set $v(r,\theta)$ as before. For a fixed $\gamma > 0$, by writing $w(t,\theta) = r^{\gamma}v(r,\theta)$ where $t = \ln r$ and recalling that $v(r,\theta)$ satisfies \eqref{polar equation}, it follows that $w(t,\theta)$ is a positive solution of
\begin{equation*}
-\partial_{t}^2 w - \Lambda_1 \partial_{t}w + \Lambda_2 w - \Delta_{\theta}w = e^{(2+\gamma - a)t}f(e^t,\theta,e^{-\gamma t} w) ~ \mbox{ in } \mathbb{R} \times \mathbb{S}^{N-1},
\end{equation*}
where
$$\Lambda_1 = N - 2 + a - 2\gamma \mbox{ and } \Lambda_2  = \gamma(N - 2 + a - \gamma).$$
If we choose $\gamma = (N-2 +a)/2$, this becomes
\begin{equation*}
-\partial_{t}^2 w + \gamma^2 w - \Delta_{\theta}w = e^{(( N+2 - a)/2) t}f(e^t,\theta,e^{-\gamma t} w) ~ \mbox{ in } \mathbb{R} \times \mathbb{S}^{N-1}.
\end{equation*}
It suffices to prove $\partial_{t}w(t,\theta) > 0$ in $\mathbb{R} \times \mathbb{S}^{N-1}$.
\medskip

{\bf Step 2.} Starting the Method of Moving Planes.

For $\lambda \in \mathbb{R}$, we set $\Sigma_{\lambda} = (-\infty,\lambda)\times \mathbb{S}^{N-1}$ and $T_{\lambda} = \partial \Sigma_{\lambda} = \{\lambda\} \times \mathbb{S}^{N-1}$. For each $t \leq \lambda$, we let $t^{\lambda} = 2\lambda - t$, which represents the reflection of $t$ across the boundary $T_{\lambda}$, and
$$ w^{\lambda}(t,\theta) = w(t^{\lambda},\theta) - w(t,\theta) \mbox{ for } (t,\theta) \in \Sigma_{\lambda}\cup T_{\lambda}.$$
By direct calculations, the comparison function $w^{\lambda}$ satisfies
\begin{align}\label{comparison inequality 1}
\partial_{t}^{2} & w^{\lambda}(t,\theta)   {}  - \gamma^2 w^{\lambda}(t,\theta)  + \Delta_{\theta}w^{\lambda}(t,\theta) \notag \\
= {} & e^{( (N+2-a)/2) t}f(e^{ t},\theta,e^{-\gamma t}w(t, \theta))   - e^{( (N+2-a)/2)  t^{\lambda}}f(e^{ t^{\lambda}},\theta,e^{-\gamma t^{\lambda}}w(t^{\lambda},\theta)) \notag \\
< {} & e^{( (N+2-a)/2)  t} [ f(e^{ t},\theta,e^{-\gamma t}w(t, \theta)) - f(e^{ t},\theta,e^{-\gamma t}w(t^{\lambda},\theta) )] ~ \mbox{ in } \Sigma_{\lambda},
\end{align}
where the last inequality follows because in $\Sigma_{\lambda}$, the subcritical growth condition \eqref{subcritical} (with $\mu = e^{2(\lambda - t)\gamma} \geq 1$) implies
$$e^{( (N+2-a)/2)  t}f(e^{t},\theta,e^{-\gamma t}w(t^{\lambda}, \theta)) < e^{( (N+2-a)/2)  t^{\lambda}}f(e^{t^{\lambda}},\theta,e^{-\gamma t^{\lambda}}w(t^{\lambda},\theta)).$$
Now we define $k(t,\theta)$ to satisfy
\begin{equation}\label{k}
k(t,\theta)w^{\lambda}(t,\theta)  =  e^{((N+2-a)/2) t}[f(e^{t},\theta,e^{-\gamma t}w(t^{\lambda}, \theta)) - f(e^{t},\theta,e^{-\gamma t}w(t,\theta) )]^{-},
\end{equation}
where $h^{-} = \min\{h, 0\}$.
Moreover, by definition, there holds
\begin{equation}\label{boundary vanish}
w^{\lambda} \equiv 0 ~\mbox{ on } T_{\lambda},
\end{equation}
and
\begin{equation}\label{mmp1}
\liminf_{t\rightarrow - \infty} w^{\lambda}(t,\theta) \geq 0 ~ \mbox{ for any fixed } \lambda \in \mathbb{R}.
\end{equation}
The comparison function satisfies
\begin{equation}\label{comparison equation}
  \begin{cases}
          \partial_{t}^2 w^{\lambda} + \Delta_{\theta} w^{\lambda} + K(t,\theta)w^{\lambda} < 0 & \text{ in } \Sigma_{\lambda}, \\
 	  w^{\lambda} = 0 & \text{ on } T_{\lambda},
        \end{cases}
\end{equation}
where $K(t,\theta) = k(t,\theta) - \gamma^2$. By the continuity of $u$, the function $e^{(-(N-2+a)/2)t}w$ is bounded in $\Sigma_{\lambda}$ uniformly in $\lambda < \bar{\lambda}$. Setting $\Sigma_{\lambda}^{-} = \{ (t,\theta) \in \Sigma_{\lambda} \, | \, w^{\lambda}(t,\theta) \leq 0\}$, for $(t,\theta) \in \Sigma_{\lambda}^{-}$ we obtain $0 \leq k(t,\theta) \leq Le^{(2+b - a)t}$ where $L$ follows from the weighted Lipschitz-type assumption on $f$. That is, we can find a positive constant $L = L(\bar{\lambda})$ such that
\[ 0 \leq k(t,\theta) \leq Le^{(2+b-a)t} \mbox{ in } \Sigma_{\lambda}^{-} \mbox{ and } \lambda < \bar{\lambda}.\]
Recalling $2 + b - a > 0$, this leads us to conclude that $\lim_{t \rightarrow -\infty} k(t,\theta) = 0$ uniformly for $(t,\theta) \in \Sigma_{\lambda}^{-}$ and $\lambda \leq \bar{\lambda}$. From this and the fact that $v$ is locally bounded, we can choose $\lambda_2 := \ln r_2$ near $-\infty$ such that for each $\lambda \leq \lambda_2$,
\begin{equation}\label{set of prop}
K(t,\theta) < -\gamma^2/2 < 0 \,\mbox{ and }\, 0 < w < 1 ~ \mbox{ in } \Sigma_{\lambda}.
\end{equation}
We apply a maximum principle argument to show
\begin{equation}\label{initiate}
w^{\lambda} \geq 0 ~ \mbox{ in } ~\Sigma_{\lambda} \, \mbox{ for all } \, \lambda \leq \lambda_2,
\end{equation}
since if otherwise, there would exist a $\lambda \leq \lambda_2$ such that $\inf_{\Sigma_{\lambda}} w^{\lambda} < 0.$ Thus by \eqref{mmp1}, $w^{\lambda}$ attains a negative minimum and due to \eqref{boundary vanish}, this minimum must be achieved away from the boundary $T_{\lambda}$. That is, there exists a point $(\bar{t},\bar{\theta}) \in \Sigma_{\lambda}$ such that
\begin{equation}\label{negativity}
w^{\lambda}(\bar{t},\bar{\theta}) = \min_{\Sigma_{\lambda}} w^{\lambda} < 0.
\end{equation}
Of course, there holds $\partial_{t}^2 w^{\lambda}(\bar{t},\bar{\theta}) + \Delta_{\theta}w^{\lambda}(\bar{t},\bar{\theta})\geq 0$. By \eqref{set of prop} and \eqref{negativity} we obtain
\begin{align*}
0 < {} & \partial_{t}^{2}w^{\lambda}(\bar{t},\bar{\theta}) + \Delta_{\theta}w^{\lambda}(\bar{t},\bar{\theta}) + K(\bar{t},\bar{\theta}) w^{\lambda}(\bar{t},\bar{\theta}),
\end{align*}
but this contradicts with \eqref{comparison equation} and therefore \eqref{initiate} holds. Clearly,
\begin{equation}\label{boundary derivative}
-\partial_{t}w^{\lambda} = 2\partial_{t} w\geq 0 \, \mbox{ on } T_{\lambda}.
\end{equation}
Moreover, the maximum principle and Hopf boundary lemma imply for any $\lambda < \lambda_2$,
\begin{equation}\label{key inequalities}
w^{\lambda} > 0 \mbox{ in } \Sigma_\lambda ~ \mbox{ and } ~ \partial_{t}w^{\lambda} < 0 \mbox{ on } T_{\lambda}
\end{equation}
since otherwise $w^{\lambda} \equiv 0$, but this cannot occur in view of \eqref{comparison equation}.

{\bf Step 3.} We show that we may continue to increase $\lambda$ so long as $w^{\lambda}$ continues to satisfy \eqref{key inequalities}. In fact, we show that we can increase $\lambda$ indefinitely. More precisely, \eqref{key inequalities} guarantees the value
$$ \lambda_0 := \sup\{ \lambda \in \mathbb{R} \, | \, w^{\mu} > 0 \, \mbox{ in } \Sigma_{\mu}, \, \partial_{t}w^{\mu} < 0 \mbox{ on } T_{\mu} ~ \mbox{ for } \mu < \lambda \}$$
exists and $\lambda_0 > -\infty$. We either have that $\lambda_0 = +\infty$ or else $\lambda_0 < +\infty$, and we claim the latter cannot happen. To this end, we assume $\lambda_0 < +\infty$.
From the growth condition \eqref{subcritical}, the fact that $w^{\lambda_0} \geq 0$ in $\Sigma_{\lambda_0}$ and \eqref{comparison inequality 1}, we get for $\lambda \leq \lambda_0$,
\begin{equation}\label{strict equation}
 \begin{cases}
          \partial_{t}^2 w^{\lambda} + \Delta_{\theta} w^{\lambda} - \gamma^2 w^{\lambda} < 0 \quad \text{ in } \Sigma_{\lambda}, \\
 	  w^{\lambda} = 0  \qquad \qquad \qquad \qquad \quad \; \text{ on } T_{\lambda}.
        \end{cases}
\end{equation}
The strong maximum principle and Hopf's boundary lemma imply that, for all $\lambda \leq \lambda_0$,
\begin{equation}\label{strong max and hopf}
w^{\lambda} > 0 ~ \mbox{ in } \Sigma_{\lambda}, ~\mbox{ and } ~ \partial_{t}w^{\lambda} < 0 ~\mbox{ on } T_{\lambda}.
\end{equation}
We shall prove there exists $\delta_0 > 0$ such that for all $0 < \delta < \delta_0$,
\begin{equation}\label{contradiction delta}
w^{\lambda_0 + \delta} \geq 0 ~ \mbox{ in } \Sigma_{\lambda_0 + \delta},
\end{equation}
thereby reaching a contradiction with the definition of $\lambda_0$.

It follows from \eqref{boundary derivative} and \eqref{strong max and hopf},that
\[ \partial_{t}w(t,\theta) > 0 ~ \mbox{ for } t \leq \lambda_0.\]
By compactness of $\mathbb{S}^{N-1}$ and continuity, we can find $\delta_1 > 0$ such that
\[ \partial_{t}w(t,\theta) > 0 ~ \mbox{ for } t \leq \lambda_0 + \delta_1,\, \theta \in \mathbb{S}^{N-1}.\]
On the other hand, we can find a sufficiently negative $\lambda_3 < \lambda_0$ such that for all $\lambda \leq \lambda_0 + \delta_1$,
\[ K(t,\theta) \leq -C < 0 ~ \mbox{ in } \Sigma_{\lambda}^{-} \cap \{t < \lambda_3\}. \]
Continuity once again guarantees we can find $0 < \delta_0 < \delta_1$ such that
\[ w^{\lambda} > 0 \mbox{ in } \Sigma_{\lambda} \cap \{\lambda_3 \leq t < \lambda\} ~ \mbox{ for all } \lambda \leq \lambda_0 + \delta_0,\]
and thus condition (II) and \eqref{k} imply
\[ k(t,\theta) = 0 \mbox{ in } \Sigma_{\lambda} \cap \{\lambda_3 \leq t < \lambda\} ~ \mbox{ for all } \lambda \leq \lambda_0 + \delta_0. \]
Thus, using similar arguments as found in Step 2, we can apply the maximum principle to conclude that \eqref{contradiction delta} holds for all  $0 < \delta < \delta_0$. Therefore, we must have $\lambda_0 = +\infty$.

Since $\lambda_0 = +\infty$, we have for each $\lambda \in \mathbb{R}$, $w^{\lambda} > 0$ in $\Sigma_{\lambda}$ and $\partial_{t}w^{\lambda} < 0$ in $T_{\lambda}$. And since $w^{\lambda}(t,\theta) = w(2\lambda - t, \theta) - w(t,\theta)$, we get
\begin{equation}\label{boundary hopf1}
\partial_{t}w^{\lambda} = -2\partial_{t} w \, \mbox{ on } T_{\lambda} ~ \mbox{ for all } \lambda \in \mathbb{R}.
\end{equation}
This leads us to $\partial_{t}w > 0$ in $\mathbb{R} \times \mathbb{S}^{N-1}$, which shows $|x|^{\gamma}u(x)$ is monotone increasing in $|x|$.  \qedhere
\end{proof}

\subsection{Some other preliminary results}
Let $1 < p < p_{S}(a,b)$, set $\alpha = \frac{2+b-a}{p-1}$ and notice that $\alpha > \gamma$. We assume the same conditions detailed in the hypotheses of Theorem \ref{thm3}.
\begin{lemma}\label{subsolution}
Fix $d > 0$ and an integer $\kappa \geq 2$. For large $R > 1$, the boundary-value problem
\begin{equation}\label{subsolution bvp}
\begin{cases}
         L_{r}[\psi(r)] \!+\! d R^{\alpha(p-1) - 2}\psi^{p}(r) \!=\! 0, \psi(r) \!>\! 0,  \mbox{ for } (\kappa - 1)R \!<\! r\! <\! (\kappa \!+\! 1)R, \\
 	  \psi((\kappa - 1)R) = \psi((\kappa + 1)R) = 0, 
        \end{cases}
\end{equation}
admits a sub-solution $\psi_{R}(r)$ with $\psi_{R}(r) = O( R^{-\alpha} )$ in $[(\kappa-1)R,(\kappa+1)R]$.
\end{lemma}

\begin{proof}
Consider the problem
\begin{equation}\label{subsolution bvp}
  \begin{cases}
         \psi''(r) + \dfrac{N-1+a}{r + \kappa}\psi'(r) + d \psi^{p}(r) = 0 & \mbox{ for } -1 < r < 1, \\
 	  \psi(-1) = \psi(1) = 0. 
        \end{cases}
\end{equation}
Indeed, for $A,K > 0$ the function $\psi_{1}(r) = K(1 - r^2)^A$ satisfies $\psi_{1}(\pm 1) = 0$ and defining $G_{1}(r) := \frac{(N-1+a)\kappa}{r+\kappa} > 0$ and
$$G_{2}(r) := 4A(A-1) - 2A(N-2+a+2A)(1-r^2) + d K^{p-1}(1-r^2)^{A(p-1) + 2},$$
 we get
\begin{align*}
\psi_{1}''(r) {} &  + \frac{N-1+a}{r + \kappa} \psi_{1}'(r) + d \psi_{1}^{p}(r) \\
= {} & K(1-r^2)^{A - 2} \big[ 2A ( (N-2+a+2A )r^2 - (N+a ) + (1-r^2)G_{1}(r) ) \\&+ d K^{p-1}(1-r^2)^{A(p-1)+2} \big] 
> {}   K(1-r^2)^{A - 2} G_{2}(r)
 \end{align*}
 in $-1 < r < 1$. If we carefully fix suitably large $A > 1$ and positive constant $K$ so that $G_{2}(r)\geq 0$, the previous calculation shows $\psi_1$ is a sub-solution of \eqref{subsolution bvp}. Now we define $\psi_R$ to be the dilated function
\begin{equation}\label{rescaled}
\psi_{R}(r) = R^{-\alpha}\psi_{1}\big( \frac{r - \kappa R}{R}\big),
\end{equation}
which satisfies $\psi_{R}((\kappa - 1)R) = \psi_{R}((\kappa + 1)R) = 0$ and
\begin{equation*}
L_{r}[\psi_R] \!=\! R^{-(\alpha+2)}\big( \psi_{1}''(\frac{r \!-\! \kappa R}{R}) \!+\! \frac{N-1+a}{ (\frac{r - \kappa R}{R}) \!+\! \kappa}  \psi_{1}'(\frac{r - \kappa R}{R})   \big) \!>\! -R^{-(\alpha + 2)}d \psi_{1}( \frac{r - \kappa R}{R} )^{p}.
\end{equation*}
Hence,
\begin{equation*}
\psi_{R}''(r) + \frac{N-1+a}{r} \psi_{R}'(r) + d R^{\alpha(p-1) -2}\psi_{R}^{p}(r) > 0 \mbox{ for } (\kappa - 1)R < r < (\kappa + 1)R.  \qedhere
\end{equation*}
\end{proof}

\begin{lemma}\label{lower bound estimate}
Let $u$ be a positive solution of \eqref{pde} in $\Omega = \mathbb{R}^N$ and let the assumptions of Theorem \ref{thm3} hold. Then for each $R > 1$
\[ u(x) \geq \min_{|y| = R} u(y) \Big(\frac{R}{|x| } \Big)^{(N-2+a)/2} ~\mbox{ in } B_{R}^{c}(0).\]
\end{lemma}
\begin{proof}
From Lemma \ref{lemma4} we get $\partial_{t} w(t,\theta) > 0$ in $\mathbb{R} \times \mathbb{S}^{N-1}$, which implies $\gamma v(r,\theta) + r \partial_{r}v(r,\theta) > 0$ in $(0,\infty)\times \mathbb{S}^{N-1}$. From the previous estimate we get $\partial_r v(r,\theta)/v(r,\theta) > -\gamma/r$ and integrating this from $R$ to $r > R$ along the radius leads to the desired conclusion.
\end{proof}


\section{Proof of Theorems \ref{thm1}--\ref{thm2} and Theorem \ref{thm4}}\label{penultimate section}

\begin{proof}[Proof of Theorem \ref{thm1}]
Let $u$ be a positive solution of \eqref{gen bvp}, and let $B_{R}(0)$ be the smallest open ball containing $\Omega$. For each $0 < \rho < R$, we define $\Sigma_{\rho} = \Omega \backslash \overline{B}_{\rho}(0)$, and for any $x \in \Sigma_{\rho}$, we define the Kelvin transform of $u$ by
$$ u_{\rho}(x) = \Big( \frac{\rho}{|x|} \Big)^{N-2+a}u(x^{\rho}), ~ x \in \Sigma_{\rho},$$
where $x^{\rho} = \rho^2 x/|x|^2$ and $u_{\rho}$ is defined in $\Sigma_{\rho}$ since $\Omega$ is assumed to be star-shaped with respect to the origin. By writing $u = v(r, \theta)$ in polar coordinates so that $v_{\rho}(r,\theta) = (\rho/r)^{N - 2 + a} v(\rho^2/r, \theta)$, direct calculations reveal that
\begin{align*}
r^{a}\Big( \partial^{2}_{r} v_{\rho}(r,\theta) +  {} & \frac{N-1+a}{r} \partial_{r}v_{\rho}(r,\theta) + \frac{1}{r^2} \Delta_{\theta}v_{\rho}(r,\theta) \Big) \notag \\
{} & + \big( \frac{\rho}{r} \big)^{N+2-a}f(\frac{\rho^2}{r}, \theta, \big(\frac{\rho}{r}\big)^{-(N-2+a)} v_{\rho}(r,\theta)) = 0.
\end{align*}
Hence, $u_{\rho}(x)$ satisfies
\begin{equation}
div(|x|^a D u_{\rho}(x)) + \Big( \frac{\rho}{|x|}\Big)^{N + 2 - a} f\Big( \frac{\rho^2 x}{|x|^2},\big(\frac{\rho}{|x|}\big)^{-(N-2+a)} u_{\rho}(x) \Big) = 0 ~ \mbox{ in } \Sigma_{\rho}.
\end{equation}
Using \eqref{supercritical} with $\mu = (|x|/\rho)^{N-2 + a} > 1$, this leads to
\begin{equation*}
div(|x|^a D u_{\rho}(x)) + f(x,u_{\rho}(x)) \leq 0 ~\mbox{ in } \Sigma_{\rho}.
\end{equation*}
Define $w_{\rho} = u_{\rho} - u$, which satisfies
\begin{equation}
 \begin{cases}
          div( |x|^{a} D w_{\rho} ) + C(x)w_{\rho} \leq 0 \quad \text{ in } \Sigma_{\rho}, \\
 	  w_{\rho} \geq 0 \qquad \qquad \qquad \qquad \qquad\, \text{ on } \partial\Sigma_{\rho},
        \end{cases}
\end{equation}
 where $C(x)$ is bounded by our assumptions on $f$. For each connected component $Z$ of $\Sigma_{\rho}$, the positivity of $u$ and since $\Omega$ is star-shaped with respect to the origin, we have $w_{\rho} > 0$ on a subset of $\partial Z\backslash \partial B_{\rho}(0)$ with positive measure. Note that we cannot assume this holds for all of $\partial \Sigma_{\rho} \backslash \partial B_{\rho}(0)$ since we do not assume $\Omega$ is strictly star-shaped. We choose small $\epsilon > 0$ so that for all $\rho \in [R - \epsilon, R)$, $meas(\Sigma_{\rho})$ is small enough so Varadhan's maximum principle for small volume domains applies and yields $w_{\rho} \geq 0$ for all such $\rho$; and the positivity of $u$ again implies that $w_{\rho} > 0$ in each $Z$ and thus on all of $\Sigma_{\rho}$.

 Next, we prove the maximal interval $(\rho_0,R)$ for which $w_{\rho} > 0$ in $\Sigma_{\rho}$ for all $\rho \in (\rho_0,R)$ is actually $(0,R)$. To this end, assume on the contrary that $\rho_0 > 0$. Note that $w_{\rho_0} > 0$ in $\Sigma_{\rho_0}$. We choose a subset $\Sigma' \subset\subset \Sigma_{\rho_0}$ and small $\epsilon \in (0, \rho_0)$ guaranteeing $meas(\Sigma_{\rho} \backslash \Sigma')$ is small enough for each $\rho \in (\rho_0 - \epsilon, \rho_0)$ so that the maximum principle for small volume domains applies to
 \begin{equation}
\begin{cases}
          div( |x|^{a} D w_{\rho} ) + C(x)w_{\rho} \leq 0 & \text{ in } \Sigma_{\rho} \backslash \Sigma', \\
 	  w_{\rho} \geq 0 & \text{ on } \partial (\Sigma_{\rho} \backslash \Sigma'), 
        \end{cases}
\end{equation}
to get $w_{\rho} \geq 0$ in $\Sigma_{\rho} \backslash \Sigma'$. Reducing $\epsilon$ if necessary, we may assume $w_{\rho} \geq \delta > 0$ in $\Sigma'$, and we conclude $w_{\rho} > 0$ in the interval $(\rho_{0} - \epsilon, R)$. This is a contradiction and thus $\rho_0 = 0$. Particularly, we have proven that for every fixed $x \in \Omega_0$,
\begin{equation}\label{comparison inequality}
\Big( \frac{\rho}{|x|}\Big)^{N-2+a} u\Big( \frac{\rho^2 x}{|x|^2} \Big) \geq u(x) ~\mbox{ for each } 0 < \rho < |x|.
\end{equation}
Sending $\rho \longrightarrow 0$ in \eqref{comparison inequality} yields an absurdity in view of the positivity of $u$.
\end{proof}

\begin{proof}[Proof of Theorem \ref{thm2}]
We shall seek radially symmetric solutions of \eqref{pde} with $\Omega = \mathbb{R}^N$ using Theorem \ref{thm1} and a simple shooting method. In particular, our strategy is to search for radially symmetric solutions, in which case \eqref{polar equation} suggests we solve the initial-value problem,
\begin{equation}\label{IVP}
\begin{cases}
\displaystyle r^a L_{r}[v(r)] + f(r,v(r)) = 0, \, v(r) > 0, \, r > 0, \medskip \\
v(0) = \beta>0, \ r^{N-1+a}v'(r) = o(1) \mbox{ at } 0.
\end{cases}
\end{equation}
The local existence and uniqueness of a positive solution $v = v(r) \in C^{2}((0,\overline{r}))\cap C([0,\overline{r}))$ for some $\overline{r} > 0$ follows by seeking fixed points of
\begin{equation}\label{integral operator}
T_{f}(v) = \beta - \int_{0}^{r}\int_{0}^{t} \Big( \frac{s^{N-1}}{t^{N-1+a}} \Big)f(s, v(s)) \,ds dt.
\end{equation}
 The latter is a consequence of our assumptions on $f$ and Banach's fixed point theorem. For convenience, we denote the solution of \eqref{IVP} by $v(r;\beta)$.

We note that $N - 1 + a > 0$ and multiplying the equation in \eqref{IVP} by $r^{N-1}$, we easily arrive at $-(r^{N-1+a} v'(r) )' > 0$ for $r > 0$.  Integrating this differential inequality leads to $v'(r;\beta) < 0$, i.e., the positive solution $v(r;\beta)$ is monotone decreasing. Now we denote by $r_0 = r_{0}(\beta) \in (0,+\infty]$ the maximal time of existence for the positive solution. Observe it is enough to show the existence of an initial shooting position $\beta > 0$ such that $r_0(\beta) = +\infty$.

To proceed, we assume the contrary that $r_{0}(\beta) < +\infty$ for each $\beta > 0$. Then the monotonicity of solutions ensures $v(r_0(\beta); \beta) = 0$. This implies that $u(x) := v(|x|;\beta) \in C^{2}(\Omega_0) \cap C(\bar{\Omega})$ is a positive (radially symmetric) solution of the Dirichlet problem \eqref{gen bvp} with $\Omega = B_{r_0}(0)$. This, however, contradicts Theorem \ref{thm1} and hence $r_0(\beta_{\ast}) = +\infty$ for some $\beta_{\ast} > 0$. Thus, $u(x) = v(|x|;\beta_{\ast})$ is a positive (radially symmetric) solution of \eqref{pde} with $\Omega = \mathbb{R}^N$.
\end{proof}

We now provide the proof of Theorem \ref{thm4}, which is mostly the same as the one for Theorems \ref{thm1}--\ref{thm2}. However, our shooting argument requires further technical modifications, since a more careful setup is needed in finding a proper initial shooting position. Namely, we need to carefully set up a degree argument to identify appropriate initial conditions, and the approach we adopt is inspired by ideas from \cite{LGZ06a}, also see \cite{LV16a,LV16b}.
\begin{proof}[Proof of Theorem \ref{thm4}]
For part (a), we suppose $U$ is a positive solution of \eqref{vector system} and we define
$$ U_{\rho}(x) = \Big( \frac{\rho}{|x|} \Big)^{N-2+a}U(x^{\rho}), ~ x \in \Sigma_{\rho},$$
where $x^{\rho}$ and $\Sigma_{\rho}$ are defined as before. Then, by adopting the same essential arguments as those in our proof of Theorem \ref{thm1} and invoking the cooperative property of $F$, we similarly conclude for each $x \in \Omega_0$,
$$ \Big( \frac{\rho}{|x|} \Big)^{N-2+a}U(x^{\rho}) \geq U(x) ~ \mbox{ for each } 0 < \rho < |x|.$$
Hence, sending $\rho \longrightarrow 0$ leads to a contradiction with the positivity of solutions.

Part (b) requires a more careful approach than the single equation case. More precisely, we must invoke a topological fixed point argument to get an appropriate initial condition that results in a global entire solution.

{\bf Step 1:} Set up an appropriate initial-value problem.

First, for each $\beta = (\beta_1, \beta_2) > 0$, we consider the initial-value problem,
\begin{equation}\label{IVP sys}
\begin{cases}
\displaystyle r^a L_{r}[v_{1}(r)] + f_{1}(r,V(r)) = 0, \, v_{1}(r) > 0, ~ r > 0, \medskip \\
\displaystyle r^a L_{r}[v_{2}(r)] + f_{2}(r,V(r)) = 0, \, v_{2}(r) > 0, ~ r > 0, \medskip \\
V(0) = \beta, \ r^{N-1+a}V'(r) = o(1) \mbox{ at } 0,
\end{cases}
\end{equation}
where recall that $V = (v_1,v_2)$. For convenience, we denote by $V(r) = V(r;\beta) \in (C^{2}((0,r_0))\cap C([0,r_0]))^2$ the unique solution of \eqref{IVP sys}, where $(0,r_0)$ is the maximal interval of existence for the positive solution. The existence, uniqueness and monotonicity of such a solution follows from standard arguments as described earlier. Then either $r_0 = r_{0}(\beta) = +\infty$, in which case $V(r;\beta)$ is indeed a global positive solution, or else $r_0 < +\infty$ and $r_0$ would be the first time at least one of the components of $V(r;\beta)$ vanishes at $r = r_0$. To complete the proof, it suffices to ultimately show the existence of a $\beta_{\ast} > 0$ such that $r_{0}(\beta_{\ast}) = +\infty$. We proceed by contradiction and assume $r_{0}(\beta) < +\infty$ for every $\beta > 0$.
\medskip

{\bf Step 2:} Introduce a ``target" map and its properties.

Recall $E =  [0,\infty)\times[0,\infty)$ and define the mapping $\varphi$ on $E$, where $\varphi(\beta) = V(r_0;\beta)$ for $r_{0} = r_{0}(\beta) < +\infty$, while we set $\varphi(\beta) = \beta$ on the boundary $\partial E$. By our definition of $\varphi$ it is easy to see $\varphi(E) \subset \partial E$.

Moreover, $\varphi: E \longrightarrow \partial E$ is indeed continuous on $E$ and we shall give a proof of this for completeness and because it illustrates the reason for imposing \eqref{source compare}. The continuity of $\varphi$ at interior points follows from basic ODE theory. So it remains to address the continuity of $\varphi$ at the boundary $\partial E$. Pick any $\beta_0 \in \partial E$. If $\beta_0 = 0$, then the positivity and monotonicity of solutions imply $|\varphi(\beta) - \varphi(\beta_0)| = |\varphi(\beta)| \leq \beta \longrightarrow \beta_0 = 0$ as $\beta \longrightarrow \beta_0$ in $E$. So let $\beta_0 \in \partial E\backslash\{0\}$ and without loss of generality, we assume $\beta_{0} = (\bar{\beta},0)$ with $\bar{\beta} > 0$. Choosing an arbitrary positive $\beta = (\beta_1,\beta_2)$ sufficiently near $\beta_0$, we may further assume $\beta_1 > \bar{\beta}/2 > 2\beta_2$. We claim that
\begin{equation}\label{prestability}
v_{2}(r;\beta) < v_{1}(r;\beta) ~ \mbox{ for } 0 \leq r < r_0 = r_{0}(\beta).
\end{equation}
If not, then there would exist a minimal $R_1 \in (0, r_0)$ such that
\begin{equation}\label{strict est}
v_{2}(r;\beta) < v_{1}(r;\beta) ~ \mbox{ in } [0,R_1) ~ \mbox{ and } ~ v_{2}(R_1;\beta) = v_{1}(R_1;\beta).
\end{equation}
In view of \eqref{source compare} and \eqref{strict est}, we see $f_{1}(r,V(r;\beta)) \leq f_{2}(r,V(r;\beta))$ for $0 \leq r < R_1$. Then
\begin{equation}\label{control}
\beta_1 - v_{1}(r;\beta) = T_{f_1}(V(r;\beta)) \leq T_{f_2}(V(r;\beta)) \leq \beta_2 - v_{2}(r;\beta)
\end{equation}
for $0 \leq r < R_1$, where $T_{f}(V)$ is defined analogously as above in \eqref{integral operator}. This implies that
\[ 0 < v_{2}(r;\beta) + \bar{\beta}/4 < v_{2}(r;\beta) + \beta_1 - \beta_2 \leq v_{1}(r;\beta) ~ \mbox{ for } 0\leq r < R_1.\]
By continuity, sending $r \longrightarrow R_1$ above leads to the absurdity $v_{2}(R_1;\beta) < v_{1}(R_1;\beta)$. Hence, the claim \eqref{prestability} holds, which further yields \eqref{control} for $0 \leq r < r_0$. This reveals that, as $\beta \longrightarrow \beta_0$, $0 \leq \beta_1 - v_{1}(r_0;\beta) \leq \beta_2 - v_{2}(r_0;\beta) \leq \beta_2 \longrightarrow 0$ and thus $ |\varphi(\beta_0) - \varphi(\beta)| = |\beta_0 - V(r_0;\beta)| \leq  |\beta - \beta_0| + |\beta - V(r_0;\beta)|\longrightarrow 0$. 

\medskip

{\bf Step 3:} We show for each real $\xi > 0$, there exists a positive $\beta_{\xi} = (\beta_{\xi,1}, \beta_{\xi,2}) \in E$ such that $\beta_{\xi, 1} + \beta_{\xi, 2} = \xi$ and $\varphi(\beta_{\xi}) = 0$.

We prove this using a simple topological degree argument.  Fix any real $\xi > 0$ and define the subsets
\[ A_\xi = \{ \beta \geq 0 \, : \, |\beta| := \beta_1 + \beta_2 = \xi \} ~ \mbox{ and } ~ B_\xi = \{ \beta \in \partial E \, : \, |\beta| \leq \xi\}. \]
Consider the homeomorphism $h:B_\xi \longrightarrow A_\xi$ defined by
\[ h(\beta) = \beta + \frac{1}{2}(\xi - |\beta|) ( 1, 1 )\]
with inverse
\[ h^{-1}(\beta) = \beta - \min\{\beta_1, \beta_2\} (1, 1).\]
The continuity of $\varphi$ implies the continuity of the composite map $h \circ \varphi : A_\xi \longrightarrow A_\xi$ on $A_\xi$ and further notice that $h \circ \varphi \equiv Identity$ on the boundary $\partial A_\xi$. Hence, for each interior point $\beta$ of $A_{\xi}$,
$$degree(h\circ \varphi, A_\xi, \beta) = degree(Identity, A_\xi, \beta) = 1$$
and the homotopy invariance of the degree implies $h\circ \varphi$ is surjective and so $\varphi$ is surjective as well. In particular, there exists a non-trivial $\beta_{\xi}$ with $|\beta_{\xi}| = \xi$ such that $\varphi(\beta_{\xi}) = 0$.
\medskip

{\bf Step 4:} Since $r_{0}(\beta_{\xi})< +\infty$ and $V(r_0;\beta_{\xi}) = \varphi(\beta_{\xi}) = 0$, it follows that $U(x) = V(|x|; \beta_{\xi})$ is a positive (radially symmetric) solution of \eqref{gen sys bvp}, where $\Omega = B_{r_0}(0)$. This contradicts with part (a) of Theorem \ref{thm4}. Therefore, we conclude there exists a positive $\beta_{\ast}$ such that $U(x) = V(|x|;\beta_{\ast})$ is an entire positive solution of \eqref{vector system} in $\Omega = \mathbb{R}^N$. This completes the proof of the theorem.
\end{proof}

\section{Proofs of Theorems \ref{thm3} and \ref{thm5}}\label{final section}
\begin{proof}[Proof of Theorem \ref{thm3}]
We proceed by contradiction and assume $u \in C^{2}(\mathbb{R}^N \backslash \{0\} ) \cap C(\mathbb{R}^N)$ is a positive solution of \eqref{pde} in $\Omega = \mathbb{R}^N$. Define $\psi_{R}$ as in Lemma \ref{subsolution} and set
$$\alpha_0 = \frac{2}{p-1} ~ \mbox{ and recall } ~ \alpha = \frac{2+b -a }{p-1}.$$
By Lemma \ref{subsolution}, Lemma \ref{lower bound estimate} and $\alpha > \gamma$, there exists a suitably large $R_0 > 1$ such that
\begin{equation*}
\psi_{R_0}(r) < u(r,\theta) \mbox{ in } \Omega_{R_0} = \big\{ (r,\theta) \, \big | \, (\kappa - 1)R_0 < r < (\kappa + 1)R_0  ,\, \theta \in \mathbb{S}^{N-1} \big \}.
\end{equation*}
By the assumptions of Theorem \ref{thm3}, we know $f(r,\theta,u) \geq d R_{0}^{b}u^p$ in $\Omega_{R_0}$ and $\alpha(p-1)-2 = b - a$. We then deduce the two inequalities
\begin{equation*}
L[u] + d R_{0}^{b - a} u^p \leq 0 \leq L_{r}[\psi_{R_0}] + d R_{0}^{b - a} \psi_{R_0}^p ~ \mbox{ in } \Omega_{R_0}.
\end{equation*}
Now, for any $0 < \delta \leq 1$, we define the rescaled functions
\begin{equation}
w_{\delta}(r) = \delta^{-\alpha_0}R_{0}^{-\alpha} \psi_{1}\big(\frac{r - \kappa R_0}{\delta R_0} \big)
\end{equation}
and set $\Omega_{\delta R_0} = \big\{ (r,\theta) \, \big | \, |r - \kappa R_0| < \delta R_0, \, \theta \in \mathbb{S}^{N-1} \big \}$.
We can easily see that $w_{1} = \psi_{R_0}$ and since $\alpha_0 (p-1) = 2$,
\begin{equation*}
  \begin{cases}
          L_{r}[w_{\delta}] + d R_{0}^{b-a}w_{\delta}^{p} \geq 0  \quad \text{ in } \Omega_{\delta R_0}, \\
 	  w_{\delta} = 0  \qquad \qquad \qquad \qquad \text{ on } \partial\Omega_{\delta R_0}.
        \end{cases}
\end{equation*}
Further noting that $w_{\delta}(\kappa R_0) \longrightarrow +\infty$ as $\delta \longrightarrow 0$, this allows us to find $0 < \bar{\delta} < 1$ and a point $(\bar{r}, \bar{\theta}) \in \Omega_{\bar{\delta} R_0}$ such that $u \geq w_{\bar{\delta}}$ in $\Omega_{\bar{\delta} R_0}$ and $u(\bar{r}, \bar{\theta}) = w_{\bar{\delta}}(\bar{r})$. As $u$ is a super-solution and $w_{\delta}$ is a sub-solution, the strong maximum principle asserts $u \equiv w_{\bar{\delta}} = 0$ on $\partial \Omega_{\bar{\delta} R_0}$. This contradicts the positivity of $u$, and this concludes the proof of the theorem.
\end{proof}

\begin{proof}[Proof of Theorem \ref{thm5}]
The proof is entirely similar to that of the proof of Theorem \ref{thm3}, so we sketch the main steps. We start by assuming $U$ is a positive solution of system \eqref{vector system} in $\Omega = \mathbb{R}^N$. The method of moving planes applies similarly to arrive at the analogous monotonicity result in the radial direction. From this monotonicity, we have for $i = 1,2$ and each $R > 1$,
\begin{equation}\label{system lower bound}
u_{i}(x) \geq \min_{|y| = R} u_{i}(y) \big(\frac{R}{|x|} \big)^{(N-2+a)/2} ~\mbox{ for } x \in B_{R}^{c}(0).
\end{equation}
Let $\Psi_{R} = (\psi_R, \psi_R)$, where $\psi_R$ is the same sub-solution derived in Lemma \ref{subsolution} with the proper choices on the parameters. By \eqref{system lower bound}, we may find a large $R_0 > 1$ such that
\[ \Psi_{R_0}(r) < U(r,\theta) ~ \mbox{ in } \Omega_{R_0} = \big\{ (r,\theta) \, \big | \, (\kappa - 1)R_0 < r < (\kappa + 1)R_0, \, \theta \in \mathbb{S}^{N-1} \big\}. \]
Now let $\alpha_0 = \max_{i=1,2} 2/(|P^i|-1)$, $\alpha = \min_{i=1,2} (2 + b - a)/(|P^i|-1)$, and by our assumptions on $F$, $ f_{i}(x,U) \geq d R_{0}^{b}U^{P^i} \mbox{ in } \Omega_{R_0}$ for each $i = 1,2$. Thus,
\[ L[u_{i}] + d R_{0}^{b-a}U^{P^i} \leq 0 \leq L_{r}[\psi_{R_0}] + d R_{0}^{b-a} \Psi_{R_0}^{P^{i}}~ \mbox{ in } \Omega_{R_0}  ~ (i = 1,2). \]
For $0 < \delta \leq 1$, define
\[ W_{\delta}(x) = \delta^{-\alpha_0} R_{0}^{-\alpha}\Big( \psi_{1} (\frac{r - \kappa R_0}{\delta R_0}) , \psi_{1} (\frac{r - \kappa R_0}{\delta R_0}) \Big) \]
and so $W_{\delta = 1} = \Psi_{R_0}$ and its components satisfy
\begin{equation}\label{contradiction bvp}
\begin{cases}
          L_{r}[w_{\delta,i}] \!+\! d R_{0}^{b-a}W_{\delta}^{P^i} \!\geq\! 0  \;\text{ in } \Omega_{\delta R_0} \!=\! \big\{ (r,\theta \, \big | \, |r- \kappa R_0| \!<\! \delta R_0, \, \theta \!\in\! \mathbb{S}^{N-1} \big\}, \\
 	  w_{i,\delta} = 0  \qquad \qquad \qquad \quad \;\text{ on } \partial\Omega_{\delta R_0}.
        \end{cases}
\end{equation}
As $w_{\delta,i}(\kappa R_0) \longrightarrow +\infty$ as $\delta \longrightarrow 0$ for $i=1,2$, this allows us to find $0 < \bar{\delta} < 1$, an index $\bar{i} \in \{1,2\}$ and a point $(\bar{r}, \bar{\theta}) \in \Omega_{\bar{\delta} R_0}$ such that $U \geq W_{\bar{\delta}}$ in $\Omega_{\bar{\delta} R_0}$ and $u_{\bar{i}}(\bar{r}, \bar{\theta}) = w_{\bar{\delta},\bar{i}}(\bar{r})$. As $U$ is a super-solution and $W_{\delta}$ a sub-solution, condition \eqref{cooperative} and the strong maximum principle imply $u_{\bar{i}} \equiv w_{\bar{\delta},\bar{i}} = 0$ on $\partial \Omega_{\bar{\delta} R_0}$. This cannot happen and thus completes the proof.
\end{proof}

\section*{Acknowledgments} This work was completed while the author was visiting Fudan University in the summer of 2019. The author would like to thank the School of Mathematical Sciences at Fudan University for their hospitality and support. The author is also partially supported by the Simons Foundation Collaboration Grants for Mathematicians 524335.


\medskip

\end{document}